\documentclass[11pt,twoside]{article}

\usepackage{longtable}
\usepackage{amssymb}

\newtheorem{Thm}{Theorem}[section]
\newtheorem{Lemma}{Lemma}[section]
\newcommand{\bfb}{\mbox{\bf B}}
\newcommand{\bfs}{\mbox{\bf S}}
\newcommand{\ds}{\displaystyle}

\begin{document}
\sloppy

\title{The Delsarte Method in the Problem of the Antipodal Contact Numbers\\
of Euclidean Spaces of High Dimensions\thanks{This work was
supported by the Russian Foundation for Basic Research, project
no. 02--01--00783, and by the Program for State Support of Leading
Scientific Schools of Russian Federation, project no.
1347.2003.1.}}

\author{D. V. Shtrom}
\date{}

\maketitle


In this paper we study the Delsarte  problem  for even functions continuous on $[-1,1]$,
nonpositive on $[-1/2,1/2],$ and representable as series with respect to the
ultraspherical polynomials $\{R_n^{\alpha,\alpha}\}_{n=0}^{\infty},$
\mbox{$\alpha=(m-3)/2,$} $m\ge 2,$ with nonnegative coefficients. The value $w^A_m$ of
the Delsarte problem gives an upper bound for the largest power of antipodal spherical
1/2-code of the space ${\mathbb{R}}^m,\ m\ge 2.$ The Delsarte problem is the problem of
infinite linear programming. In this paper the value $w^A_m$ is found for $3\le m \le
161$ with a few gaps. The results are summarized in a table. We give a detailed proof for
$m=43$ and point out a scheme of a proof for other cases.

\section{Introduction. Setting of the problem} In this
paper we study an extremal problem on a class of functions continuous on a closed
interval, representable by series with respect to the ultraspherical polynomials with
restrictions on values of the functions and coefficients of the representations. The
problem occurs  in applying of the Delsarte scheme to the problem on the largest power of
antipodal spherical $s$-codes of Euclidean spaces $\mathbb{R}^m.$

For the first time this scheme arose in the investigations of Delsarte \cite{DEL1},
\cite{DEL2}  of the bounds of packings in certain metric spaces. The Delsarte scheme was
developed and successfully used   in the works of G.A.~Kabatianskii and V.I.~Levenshtein
\cite{Kabatiansky and Levenshtein}, A.~Odlyzko and N.~Sloane \cite{Odlyzko and Sloane},
V.I.~Levenshtein \cite{Levenshtein}, \cite{Lev-83}, \cite{Lev-92}, V.M.~Sidel'nikov
\cite{Sidelnikov}, V.A.Yudin \cite{UdinAnd}, P.G.~Boyvalenkov \cite{BoivalenkovNE},
P.G.~Boyvalenkov, D.P.~Danev, and S.P.~Bumova \cite{BoivDan1996}, V.V.Arestov and
A.G.~Babenko \cite{Arestov and Babenko}, \cite{Arestov and Babenko - diktatory}, and
others in connection with the investigation of optimal (in one or another sense)
arrangement of points in metric spaces and, in particular, the investigation of the
contact numbers of Euclidean spaces $\mathbb{R}^m$. The Delsarte scheme leads to a
problem of infinite linear programming (we will call it the Delsarte problem). The value
of the problem gives an upper bound for the original problem. Our interest is the
Delsarte problem connected with spherical $s$-codes in $\mathbb{R}^m.$ First in this
theme, V.M.Sidel'nikov used in his work a quadrature formula of the Gauss--Markov type
where only values of a function were equipped.  He applied the formula to prove an
extremality of the Levenshtein polynomials on a class of polynomials of the same or a
smaller degree. In the paper \cite{BoivDan1996} by P.G.Boyvalenkov, D.P.Danev, S.P.Bumova
an extremality criterion for the polynomials of V.I.Levenshtein in the problem on the
largest power  of a spherical $s$-code was proved. In the paper \cite{BoivDan1998} by
P.G.Boyvalenkov and  D.P.Danev, the criterion was generalized. This made possible the use
of the criterion to establish values $s$ and $m$ for which the Levenshtein polynomials
are the extremal polynomials in the problem on the largest power of an antipodal $s$-code
(the antipodal contact number). Using this criterion, P.G.Boyvalenkov and D.P.Danev
established the non-extremality of the Levenshtein polynomials for some $s$ and $m.$ But
they failed to establish the extremality of the Levenshtein polynomials. V.V.Arestov and
A.G.Babenko \cite{Arestov and Babenko} constructed essentially a quadrature formula where
not only values of functions, but also their Fourier--Jackobi coefficients were equipped.
In this way they succeeded, in particular, in solving of the Delsarte problem for $s=1/2,
m=4$. In the paper \cite{Arestov Babenko and Dekalova} the Delsarte problem was solved
for $s=1/3, \, m=4,\,5,\,6$ by a similar method. The author \cite{Shtrom1} applied this
method to find a solution of the Delsarte problem for $s=1/2$ and the following $m:$
  $$ 5\leq m\leq146,\ 148\leq m\leq156,\ m=161, \ (m\ne 8, 24).$$

Let $\mathbb{R}^m,\, m\ge2$,\, be the real Euclidean space with the standard inner
product  $xy\,=\,x_1y_1\,+\,x_2y_2\,+\,\ldots\,+\,x_my_m,$\
\mbox{$x=(x_1,\,x_2,\,\ldots,x_m),$} $y=(y_1,x_2,\ldots,y_m),$ and the norm
$|x|=\sqrt{xx},\ x,y\in\mathbb{R}^m$.  Let $\bfb^m(y)=\{x\in\mathbb{R}^m: |x-y|\le 1\}$
be the ball of the unit radius with center at the point $y\in\mathbb{R}^m$. Let $\tau_m$
denote the maximal number $\tau$ of nonoverlapping balls $\bfb^m(y^{(1)}),$
$\bfb^m(y^{(2)}), \ \ldots,$ $\bfb^m(y^{(\tau)})$ of the unit radius: ($|y^{(i)}|=2,$
$i=1,2,\ldots,\tau;$ \ $|y^{(i)}-y^{(j)}|\ge 2,\  i\not=j$) touching the central ball
$\bfb^m(0)$ of the unit radius. The quantity $\tau_m$ is called the {\it contact number}
of the space $\mathbb{R}^m.$

At present, the exact values of $\tau_m$ are known (see \cite[Table 1.5 and
Chapter\nobreakspace 1]{KONSL}) only for $m=2,\, 3,\, 8,\, 24$; namely  $\tau_2=6,$
$\tau_3=12,$ $\tau_8=240,$ \mbox{$\tau_{24}=196560.$} In the general case two-sided
estimates of this value exist. Under additional restriction of symmetry with respect to
origin on a location of balls, we obtain the problem of non-lesser interest about finding
of {\it antipodal contact number} $\tau^A_m.$ This problem is connected (see, for
instance, \cite{UdinAnd}) with the classical problem on the largest number of integer
points on an ellipsoid. About these problems it is known a little bit more. In
particular, it is known that $\tau^A_2=6,$ $\tau^A_3=12,$ $\tau^A_4=24,$ $\tau^A_5=40,$
$\tau^A_6=72,$ $\tau^A_7=126,$ $\tau^A_8=240,$ $\tau^A_{24}=196560.$ Moreover, 24 is the
largest dimension where these values are calculated exactly. For $m>24$ only two-sided
estimates of these quantities are known. But the known estimates do not show, for
instance, the order of growth neither of $\tau_m$ nor of $\tau^A_m$  as $m\to\infty$.

Specific arrangements of balls give lower bounds for $\tau^A_m$ and $\tau_m.$ A
nonconstructive method to estimate $\tau^A_m$ from below also exists. In this paper we
will not discuss lower bounds for $\tau^A_m.$ The information and the bibliography
concerning this theme can be found in the already cited monograph \cite{KONSL}. The
Delsarte approach give an efficient upper bound for $\tau^A_m$. We will state the method
for a more general case  for the problem on the largest power of spherical $s$-codes of
Euclidean spaces $\mathbb{R}^m.$

Let $\bfs^{m-1}=\{x\in\mathbb{R}^m: |x|=1\}$ be the unit sphere of the space
$\mathbb{R}^m.$ We basically follow the notation of \cite{Delsarte Goethalts Seidel}.
Suppose a set $W\subset\bfs^{m-1}$ contains at least two points.  All possible values of
inner product of different vectors $x, y$ from $W$ form the set which we denote by ${\cal
A}(W)$, i.e.,
 $$
 {\cal A}(W)=\{xy:\quad x\not=y, \ x, y\in W\}.
 $$
The set $W\subset\bfs^{m-1}$ with the property ${\cal A}(W)\subset [-1,s]$ is called the
{\it spherical $s$-code.} The set $W\subset\bfs^{m-1}$ is called the {\it
centrosymmetrical set,} if for any point $x$ from this set, the point $-x$ also belongs to
the set. The set $W\subset\bfs^{m-1}$ which is the centrosymmetrical set and the
spherical $s$-code simultaneously, is called {\it the antipodal spherical $s$-code.}

Let ${\cal W}^A_m(s)$ denote the set of all antipodal spherical \mbox{$s$-codes} from
$\bfs^{m-1}$. Let $N^A_m=N^A_m(s)$ denote the largest possible power of antipodal
spherical $s$-code, i.e., suppose that
 \begin{equation}\label{v1}
 N^A_m(s)=\max\,\{{\rm{card}}(W): \  W\in{\cal W}^A_m(s)\}.
 \end{equation}
It is readily seen that $N^A_m(1/2)$ coincides with
 $\tau^A_{m}:$
 \begin{equation}\label{v2}
 \tau^A_{m}=N^A_m(1/2),\quad m\ge 2.
 \end{equation}

A trivial connection between $\tau_m$ and $\tau^A_m$ follows from the fact that an
antipodal spherical $s$-code is a spherical $s$-code. Namely,
$$\tau^A_m\le\tau_m.$$

Let $R_k=R_k^{\alpha,\alpha}, \ k=0,1,2,\ldots\,,$ be the system
of ultraspherical polynomials (Gegenbauer polynomials)
orthogonal on the closed interval $[-1,1]$ with weight
 $(1-t^2)^\alpha,$ $\alpha=(m-3)/2,$ and normalized by the condition $ R_k(1)=1.$
We denote by ${\cal F}^A_m={\cal F}^A_m(s),\ s\in[0,1),$ the set
of functions $f$ continuous on $[-1,1]$ with the following
properties:

 (1) the function $f$ can be represented as the series
 \begin{equation}\label{v3}
 f(t)=\sum_{k=0}^\infty f_{2k} R_{2k}(t),
 \end{equation}
whose coefficients satisfy the conditions
 \begin{equation}\label{v4}
 f_0>0, \quad  f_{2k}\ge 0, \  k=1,2,\ldots, \quad  f(1)=\sum_{k=0}^\infty f_{2k}<\infty;
 \end{equation}

 (2) the function $f$ is nonpositive on $[-s,s]:$
 \begin{equation}\label{v6}
 f(t)\le 0,\quad t\in [-s,s].
 \end{equation}
On this set of functions we consider the Delsarte problem of
evaluating
 \begin{equation}\label{v7}
w^A_m(s)=\inf\left\{\frac{f(1)+f(-1)}{f_0}: \ f\in{\cal
F}^A_m({\rm s})\right\}=\inf\left\{\frac{2\, f(1)}{f_0}: \
f\in{\cal F}^A_m({\rm s})\right\}.
 \end{equation}
Let agree to call this quantity the Delsarte constant (function).

In the paper \cite{Arestov and Babenko} by V.V.Arestov and A.G.Babenko, it was proved
finiteness of the expansion (\ref{v3}) of the function extremal in the problem
(\ref{v7}), i.e., it was proved that the function providing the minimum in the problem
(\ref{v7}) is a polynomial. At the same time, using Corollary 2.2 from \cite{Arestov and
Babenko}, point-wise estimates for the polynomials  $R_{2k}(t)$ on $[-s,s]$ (one can use,
for instance, Lemma 2.1 from \cite{Shtrom1}), and also an universal upper estimate by
V.I.Levenshtein for $w_m^A(s)$, one can find an efficient upper  estimate of an degree of
an extremal polynomial.

The following statement is  contained in the paper \cite{UdinAnd}. This statement gives
an upper estimate of $N^A_m(s)$,  and so of $\tau^A_{m}$ (see (\ref{v2})), by (\ref{v7}).

 \ \

 \noindent {\bf Theorem A.} {\it For any $s\in[0,1),$ $m=2,3,\ldots$ we have}
 \begin{equation}\label{v8}
 N^A_m(s)\le w^A_m(s).
 \end{equation}

Since $N^A_m$ is an even number, the theorem involves
 \begin{equation}\label{v9}
 N^A_m(s)\le [w^A_m(s)]_2,\quad s\in[0,1),\ \ m=2,3,\ldots\,,
 \end{equation}
where $[t]_2$ is the largest even number not exceeding $t$.

The problem on the largest power of antipodal spherical codes was considered in the
papers \cite{Lev-83}, \cite{Lev-92}, \cite{Boivalenkov3}, \cite{BoivalenkovNE}. In the
paper by V.I.Levenshtein a good upper estimate of $w_m^A$ was obtained. Bellow in the present
paper it will be shown that for rather large number of $m$ this estimate can not be improved by the
Delsarte method. In the papers \cite{BoivDan1998} and \cite{Boivalenkov3} necessary and
sufficient conditions for extremality of the Levenshtein polynomials in the Delsarte
problem were established.  In the paper \cite{Boivalenkov3} with use of these conditions
(criterion), a number of $s$ and $m$ for which the Levenshtein polynomial is not the
extremal polynomial in the Delsarte problem were found. With use of the fact, in the
paper \cite{Boivalenkov3} the estimate by V.I.Levenshtein was a little bit improved for
these $s$ and $m$. But neither in this paper nor in more earlier papers
extremality of the Levenshtein polynomials  (corresponding sufficient conditions) was not
under consideration with the exception of those cases when an upper estimate coincides
with the known lower estimate. In these cases extremality of the Levenshtein polynomials
is obtained automatically (for example, for $s=1/2$ and $m=2, 4, 6, 7, 8, 24$).

Our interest is $s=1/2.$ Let $$w^A_m=w^A_m(1/2).$$

The following is done in this paper.

(1)  Exact values of the Delsarte constant $w^A_m$ are evaluated for many $m$; namely,
for all natural \mbox{$3\leq m\leq99,$} $104\leq m\leq122,$ \mbox{$125\leq m\leq134,$}
$136\leq m\leq145,$ $147\leq m\leq156,$ $m=161$ ($m\ne 4, 6, 7, 8, 24$). Corresponding
extremal polynomials are found. Notice that polynomials of the forms 1 and 2 (see Sect.
4) are (to within a change of variable) the polynomials obtained in the paper by
V.I.Levenshtein \cite{Lev-92}. For all these $m$ it was, practically, needed to verify
sufficient conditions for extremality of the Levenshtein polynomials in the Delasarte problem.
And it was done.
The corresponding results are given in Table 1. In addition, extremal polynomials
different from the Levenshtein polynomials were found. These are polynomials of the forms
3 and 4 (see Sect. 4). The corresponding results are given in Table 2.

(2) We give a detailed proof for $m=43$. The other results are given in a table without
detailed proof because the corresponding computations occupy too much room. But the table
allows one to carry out the computations for any $m$ indicated in it similarly to the
case $m=43.$

By now, a few $s$ and $m$ are known for which $w^A_m(s)$ is an integer and coincides with
$N^A_m(s),$ for example, for $s=1/2$ and $m$ equal to $2, 4, 6, 7, 8, 24.$ But whether
the opposite is true? And $N^A_m(s)$ coincides at least with $[w^A_m(s)]_2$ whether or
not? Using the results of the paper by P.G.Boyvalenkov \cite{BoivalenkovNE} and this
paper, one can respond to the question as follows: given quantities are not equal in
general. For example, $\tau^A_5=N^A_5(1/2)=40$ but $w^A_5=42.$ It is also known
\cite{BoivalenkovNE} that $\tau^A_{10}\le 548,$ $\tau^A_{14}\le 2938,$ while
$w^A_{10}=550$ and $w^A_{14}=2940.$ Thus, we have $\tau_m^A<[w_m^A]_2$ for $m= 5, 10,
14.$

\section{\bf A method of investigation for the Delsarte problem.} From now on we follow
in general the scheme of reasoning from the papers \cite{Arestov and Babenko},
\cite{Arestov and Babenko - diktatory}, \cite{Arestov Babenko and Dekalova},
\cite{Shtrom1}. The present paper is close to the paper \cite{Shtrom1} with respect to
methods of investigation and nature of results.

Suppose $\ell=\ell_1$ is the space of summable sequences $x=\{x_{2k}\}_{k=1}^{\infty}$ of
real numbers, and $C[-1,1]$ is the space of functions continuous on $[-1,1]$. Let $A$ the
linear operator from $\ell=\ell_1$ to $C[-1,1]$ defined by
$$(Ax)(t)= \sum_{k=1}^{\infty}x_{2k} R_{2k}(t),\ \ t\in[-1,1],\
x=\{x_{2k}\}_{k=1}^{\infty}\in \ell.$$ Suppose
   \begin{equation}\label{p3_0}
 u^A_m(s)=\inf \left\{\sum_{k=1}^{\infty}x_{2k} :\ x_{2k}\geq 0;\ 1+(Ax)(t)\leq 0,\  t\in
 [-s,s]\right\}.
 \end{equation}
In the problem (\ref{v7}) we can restrict our attention to functions $f\in {\cal F}^A_m,$
with $f_0=1.$ For such a function let $x=\{f_{2k}\}_{k=1}^\infty;$ we have
$f(t)=1+(Ax)(t)$ and $f(1)=1+\sum\limits_{k=1}^{\infty}f_{2k}.$ Hence, we can deduce that
(\ref{v7}) is connected with (\ref{p3_0}) by
 \begin{equation}\label{p3_1}
 w^A_m(s)=2+2u^A_m(s).
 \end{equation}
The problems (\ref{v7}),(\ref{p3_0}) are problems of infinite linear programming (see,
for instance, the monograph \cite{golnsht}). In \cite{Arestov and
 Babenko}, with the help of these
considerations the existence of solutions (extremal functions) of a problem more general
than (\ref{v7})  and of the corresponding dual problem  was proved.

We denote by $\Phi_m^A$ the set of even functions $f\in C[-1,1]$ representable as the
series with respect to the ultraspherical polynomials $R_{2k}=R_{2k}^{\alpha,\alpha},\
\alpha=(m-3)/2,$ with an (absolutely) summable sequence of real coefficients:
 \begin{equation}\label{h0A}
 \Phi^A_m=\Big\{f\in C[-1,1]: \quad f(t)=\sum_{k=0}^\infty f_{2k} R_{2k}(t), \ \
 \sum_{k=0}^\infty |f_{2k}|<\infty\Big\}.
 \end{equation}
Evidently, ${\cal F}^A_m\subset \Phi^A_m.$

In what follows, we consider only the case $s=1/2.$ To solve the problem (\ref{v7}), a
quadrature formula (specific for each $m$) on the class of functions $f\in\Phi^A_m$ was
constructed. This formula contains not only the values of the function, but also the
Fourier coefficients $f_{2\nu}$ of the representation $ f(t)=\sum_{k=0}^\infty f_{2k}
R_{2k}(t)$ as the series with respect to the polynomials $R_{2k}$. To be exact, the
formula has the form
 \begin{equation}\label{h1}
 f_0=\frac{1}{\vartheta(\alpha)}\int _{-1}^1 f(t)(1-t^2)^{\alpha}dt = L(f)-\sum_{\nu\ge 1} L(R_{2\nu})f_{2\nu},
 \end{equation}
where
 $$
 \vartheta(\alpha)=\int_{-1}^1 (1-t^2)^{\alpha}dt, \quad \alpha=(m-3)/2,
 $$
and the functional $L$ is given by
 \begin{equation}\label{h2}
  L(f)=\lambda(1)f(1)+\lambda\left(\frac{1}{2}\right)f\left(\frac{1}{2}\right)
 +\lambda(t_1)f(t_1)+\ldots+\lambda(t_{k(m)})f(t_{k(m)}),
 \end{equation}
where the nodes \ $t_1, t_2, \ldots, t_{k(m)}$ are from $[0,1/2)$, the coefficients
$\lambda\left(\ds\frac{1}{2}\right),$
$
 \lambda(t_1),$ $
 \lambda(t_2),\ldots,\lambda(t_{k(m)})$ and the values $L(R_{2\nu}),\ \nu \ge 1,$
of the functional $L$ are nonnegative for the polynomials $R_{2\nu}$. Simultaneously with
the quadrature formula (\ref{h1})  we construct a polynomial $f^*\in{\cal F}^A_m$ such
that the following equality holds:
 \begin{equation}\label{h3}
 f^*_0=\frac{1}{\vartheta(\alpha)}\int _{-1}^1 f^*(t)(1-t^2)^{\alpha}dt =
 \lambda(1)f^*(1).
 \end{equation}
Under these conditions, the quantity $2/\lambda(1)$ is the value of the problem
(\ref{v7}) (for $s=1/2$), i.e., we have
 \begin{equation}\label{h4}
 w^A_m=\frac{2}{\lambda(1)}.
 \end{equation}

Indeed, by (\ref{h1}) and the nonnegativity of the coefficients of the quadrature
formula, for any function $f\in {\cal F}^A_m$ we have
$$f_0=L(f)-\sum\limits_{\nu\ge 1} L(R_{2\nu})f_{2\nu}\leq L(f)\leq
\lambda(1)f(1).$$ Thus, for any function $f\in {\cal F}^A_m$ we have
\begin{equation}\label{h5}
\frac{f(1)}{f_0}\geq \frac{1}{\lambda(1)}
\end{equation}
and so
\begin{equation}\label{h6}
w^A_m\ge \frac{2}{\lambda(1)}.
\end{equation}
For the polynomial $f^*$ inequality (\ref{h5}) turns into equality. Consequently,
inequality (\ref{h6}) is in fact an equality, i.e., (\ref{h4}) holds. In addition, the
function $f^*$ has the property $w^A_m=2\, f^*(1)/f^*_0,$ i.e., this polynomial is a
solution (an extremal function) of  the problem (\ref{v7}). The functional $L$ is defined
by the measure which is a solution of the dual problem.

The condition (\ref{h3}) imposes rather severe restrictions on the
function
\begin{equation}\label{hf8} f^*(t)=\sum_{k=0}^{n(m)} f^*_{2k} R_{2k}(t),
\end{equation}
and on the quadrature formula (\ref{h2}). Namely, the following conditions must be
satisfied:

(a) all nodes $t_\nu$ of the functional $L$ (except $t=1$) belong to $[0,1/2]$ and they
are zeros of the function $f^*.$ Moreover, every zero from $[0,1/2)$ is at least a double
zero;

(b) for $k\ge 1$ the coefficients $f^*_{2k}$ of the representation (\ref{hf8}) of the
function  $f^*$ are connected with the values $L(R_{2k})$ of the functional (\ref{h2})
for the polynomials $R_{2k}$ by the relation $f^*_{2k}\, L(R_{2k})=0;$

(c) the sum of weights of the functional $L$ is equal to 1, i.e., $$L(R_0)=
\lambda(1)+\lambda\left(\frac{1}{2}\right)
 +\lambda(t_1)+\lambda(t_2)+\ldots+\lambda(t_{k(m)})=1.$$

Besides, in order to the functional $L$ be nonnegative and the polynomial $f^*$ belong to
the class, the following conditions must be satisfied:

(d) the weights $\lambda (t_\nu )$ of the functional $L$ are nonnegative and $\lambda
(1)>0;$

(e)  for all polynomials $R_{2k},\, k\ge 1,$ the functional $L$ is nonnegative:
$L(R_{2k})\ge 0,\, k\ge 1;$

(f) the polynomial $f^*$ belongs to the class ${\cal F}^A_m,$ i.e., its coefficients
$f^*_{2k},\ k\ge 1,$ are nonnegative, $f^*_0>0,$ and the condition $f^*(t)\le 0, \ t\in
[-1/2,1/2],$ holds.

In all cases when the problem  (\ref{v7}) is solved exactly, $t=1/2$ is a zero of the
polynomial, moreover, its multiplicity is 1. In order to construct the polynomial $f^*$,
it is important to have a priori information about the structure of the polynomial: the
degree,  the numbers of  the vanishing coefficients from the representation (\ref{hf8}),
the number of (multiple) zeros of the polynomial on $(-1/2,1/2),$ and whether or not  the
point $t=-1$ is a zero of the polynomial. Information about the form of the extremal
function must be cleared up separately. Most often we actually guess it, using some
ideas, for example, a preliminary numerical experiment. If we have the above information,
the extremal polynomial $f^*$ and the functional $L$ are constructed as follows.

The conditions (a)--(c) give the system of (nonlinear) equations with respect to the
nodes $\{t_i\}$ of the functional (\ref{h2}) from $[0,1/2]$ that are zeros of the
function $f^*$, the coefficients $\{f^*_{2k}\}$ of the function (polynomial) (\ref{hf8}),
and the weights $\{\lambda(t_i)\}$ of the functional (\ref{h2}). Usually, a solution of
the set of equations is not unique.  We have to choose a solution that satisfies the
above-listed conditions (d)--(f).

 \section{The forms of extremal functions.}

In this paper, the value  of  $w^A_m=w^A_m(1/2)$ will be found for the following $m:\ $
$3\leq m\leq99,$ $104\leq m\leq122,$ $125\leq m\leq134,$ $136\leq m\leq145,$ $147\leq
m\leq156,\ m=161$ ($m\ne 4, 6, 7, 8, 24$). In all these cases, the extremal function is a
polynomial. The structure of the polynomial essentially depends on $m.$ The extremal
polynomials are of one of the following 4 forms. In formulas below $a_i$ denote (double)
zeros $t_i$ of an extremal polynomial (nodes of the functional $L$) belonging to $(0,
1/2).$ In these formulas, and in what follows, $s=1/2.$

 1. $f(t)=(t^2-s^2)\cdot\prod\limits_{i=1}^K(t^2-a^2_i)^2,$

 2. $f(t)=(t^2-s^2)\cdot\prod\limits_{i=1}^K(t^2-a^2_i)^2\cdot t^2,$

 3. $f(t)=(t^2-s^2)\cdot\prod\limits_{i=1}^K(t^2-a^2_i)^2\cdot(t^4+q\,t^2+r), \quad f_{4K+2}=f_{4K+4}=0,$

 4. $f(t)=(t^2-s^2)\cdot\prod\limits_{i=1}^K(t^2-a^2_i)^2\cdot(t^4+q\,t^2+r)\cdot t^2, \quad f_{4K+4}=f_{4K+6}=0,$

\ \

The results of the work are summarized in the Tables 1 and 2 below. The first column
gives the space dimension, the third column gives the number $K$ of double zeros of the
extremal function on $(0,1/2)$, the fourth column gives one of 4 forms (exactly, the
number of the form) of the extremal function, and the fifth column gives the degree of
the obtained polynomial. The second column contains the solution of the Delsarte problem
$w^A_m=w^A_m(1/2)$  that we have found.

The form of the extremal function and the values $K$, $m$ is the very information
according to that we form  the simultaneous equations starting from the conditions
(a)--(c) of Section~2. In all cases considered in the paper, the number of equations
coincides with the number of variables. The equations are nonlinear. In this connection,
the solution is not unique. We have to choose the solution that gives a function $f^*$
and a functional $L$ satisfying the conditions (d)--(f). We have done this for all the
above-mentioned $m.$ But we cannot give full proof for each case because the
corresponding computations  occupy too much room. A detailed proof will be given only for
$m=43.$  Constructing of the simultaneous equations for other $m$, their analysis and
proof of extremality of obtained solutions are realized analogously.

All analytic  and numerical computations were made by using the Maple package of analytic
computations.

\section{\bf The table of the results.} \nopagebreak

In the Table 1 the values $w^A_m$ are given for the case of an extremal polynomial having
the form 1 or 2 (i.e., the Levenshtein polynomial). The indicated values were known
earlier as estimates of the largest power of an antipodal spherical $1/2$-code (see
\cite[p.8]{Lev-92}). In the present paper we assert that the mentioned values are a
solution of the problem (\ref{v7}) (for $s=1/2$).

\begin{center}
Table 1
\end{center}
\begin{longtable}{|@{}c|@{}c|@{}c|@{}c|@{}c|}
 \hline
$\phantom{n}m\phantom{n}$& $w^A_m$ & $\,K\,$ &\mbox{ Form}&\mbox{ Degree}
\\ \hline
\endfirsthead
\caption[Table 1]{(continued)}\\ \hline $\phantom{n}m\phantom{n}$& $w_m$ & $\,K\,$
&\mbox{ Form}&\mbox{ Degree}
\\
\hline
\endhead
 \hline 4 & $24$ &0& 2 & 4
 \\ \hline 5 & $42$ &0& 2 & 4
 \\ \hline 6 & $72$ &0& 2 & 4
 \\ \hline 7 & $126$ &0& 2 & 4
 \\ \hline 8 & $240$ &0& 2 & 4
 \\ \hline 9 & $366{\frac {12}{73}}$ &1& 1 & 6
 \\ \hline 10 & 550 &1& 1 & 6
 \\ \hline 11 & $820{\frac {16}{23}}$ &1& 1 & 6
 \\ \hline 12 & $1228{\frac {1}{2}}$ &1& 1 & 6
 \\ \hline 13 & $1867{\frac {17}{19}}$ &1& 1 & 6
 \\ \hline 14 & 2940 &1& 1 & 6
 \\ \hline 15 & $4962{\frac {6}{37}}$ &1& 1 & 6
 \\ \hline 16 & 8160 &1& 2 & 8
 \\ \hline 17 & $11478{\frac {12}{13}}$ &1& 2 & 8
 \\ \hline 18 & $16122{\frac {6}{7}}$ &1& 2 & 8
 \\ \hline 19 & 22724 &1& 2 & 8
 \\ \hline 20 & 32340 &1& 2 & 8
 \\ \hline 21 & $46879{\frac {7}{17}}$ &1& 2 & 8
 \\ \hline 22 & $70165{\frac {1}{3}}$ &1& 2 & 8
 \\ \hline 23 & $111126{\frac {6}{19}}$ &1& 2 & 8
 \\ \hline 24 & 196560 &1& 2 & 8
 \\ \hline 25 & $267628{\frac {622}{1451}}$ &2& 1 & 10
 \\ \hline 26 & 364182 &2& 1 & 10
 \\ \hline 27 & $497035{\frac {7}{281}}$ &2& 1 & 10
 \\ \hline 28 & 683240 &2& 1 & 10
 \\ \hline 29 & $951235{\frac {31}{251}}$ &2& 1 & 10
 \\ \hline 30 & $1352089{\frac {1}{71}}$ &2& 1 & 10
 \\ \hline 31 & $1987341{\frac {15}{197}}$ &2& 1 & 10
 \\ \hline 32 & $3091334{\frac {2}{5}}$ &2& 1 & 10
 \\ \hline 34 & 7314012 &2& 2 & 12
 \\ \hline 35 & $9768755{\frac {475}{733}}$ &2& 2 & 12
 \\ \hline 36 & 13090896 &2& 2 & 12
 \\ \hline 37 & $17663588{\frac {56}{135}}$ &2& 2 & 12
 \\ \hline 38 & $24107066{\frac {2}{59}}$ &2& 2 & 12
 \\ \hline 39 & $33491675{\frac {1}{49}}$ &2& 2 & 12
 \\ \hline 40 & $47830565{\frac {5}{47}}$ &2& 2 & 12
 \\ \hline 41 & $71400259{\frac {1}{179}}$ &2& 2 & 12
 \\ \hline 42 & $115143336$ &2& 2 & 12
 \\ \hline 44 & 238814520 &3& 1 & 14
 \\ \hline 45 & $315542890{\frac {610}{1271}}$ &3& 1 & 14
 \\ \hline 46 & $419023338{\frac {46}{53}}$ &3& 1 & 14
 \\ \hline 47 & $561215167{\frac {4037}{5461}}$ &3& 1 & 14
 \\ \hline 48 & $761656254{\frac {6}{11}}$ &3& 1 & 14
 \\ \hline 49 & $1054475186{\frac {8}{17}}$ &3& 1 & 14
 \\ \hline 50 & $1504942258{\frac {16}{23}}$ &3& 1 & 14
 \\ \hline 51 & $2255135531{\frac {583}{1363}}$ &3& 1 & 14
 \\ \hline 55 & $9437927703{\frac {133}{229}}$ &3& 2 & 16
 \\ \hline 56 & $12438770728{\frac {8}{43}}$ &3& 2 & 16
 \\ \hline 57 & $16538544622{\frac {31216}{39785}}$ &3& 2 & 16
 \\ \hline 58 & $22282754713{\frac {1}{23}}$ &3& 2 & 16
 \\ \hline 59 & $30616778153{\frac {137891}{168533}}$ &3& 2 & 16
 \\ \hline 60 & 43329012480 &3& 2 & 16
 \\ \hline 61 & $64250386884{\frac {27108}{41039}}$ &3& 2 & 16
 \\ \hline 66 & $350194168928{\frac {6406}{8155}}$ &4& 1 & 18
 \\ \hline 67 & $461936954642{\frac {6152}{63809}}$ &4& 1 & 18
 \\ \hline 68 & $616741349591{\frac {433}{1237}}$ &4& 1 & 18
 \\ \hline 69 & $838108246881{\frac {6081}{9829}}$ &4& 1 & 18
 \\ \hline 70 & 1169132164200 &4& 1 & 18
 \\ \hline 71 & $1698060388955{\frac {28265}{298283}}$ &4& 1 & 18
 \\ \hline 77 & $12411939766938{\frac {246}{511}}$ &4& 2 & 20
 \\ \hline 78 & $16405689281448{\frac {24}{37}}$ &4& 2 & 20
 \\ \hline 79 & $22010329107332{\frac {140248}{248905}}$ &4& 2 & 20
 \\ \hline 80 & $30177990957237{\frac {17}{19}}$ &4& 2 & 20
 \\ \hline 81 & $42744922122472{\frac {42898}{53507}}$ &4& 2 & 20
 \\ \hline 82 & $63758049542988{\frac {64}{65}}$ &4& 2 & 20
 \\ \hline 89 & $561945080967167{\frac {37795}{1048459}}$ &5& 1 & 22
 \\ \hline 90 & $757112348026609{\frac {2267}{6637}}$ &5& 1 & 22
 \\ \hline 91 & $1045988435188294{\frac {5400226}{9600141}}$ &5& 1 & 22
 \\ \hline 92 & $1501230731871410{\frac {2}{257}}$ &5& 1 & 22
 \\ \hline 112 & $816451378740698787{\frac {5}{9}}$ &6& 1 & 26
 \\ \hline 113 & $1136254535300176728{\frac {151306698}{223651759}}$ &6& 1 & 26
 \\ \hline 114 & \phantom{1}$1136254535300176728{\frac {151306698}{223651759}}$\phantom{1} &6& 1 & 26
 \\
\hline
\end{longtable}

In the table 2 the values $w^A_m$ are given for the case when an extremal polynomial in
the problem (\ref{v7}) (for $s=1/2$) is a polynomial of the form 3 or 4. These values can
be used as new estimates of the largest power of spherical $1/2$-codes in
${\mathbb{R}}^m$ which are unimprovable by the Delsarte method.

\newpage

\begin{center}
Table 2
\end{center}
\begin{longtable}{|@{}c|@{}c|@{}c|@{}c|@{}c|}
 \hline
$\phantom{n}m\phantom{n}$& $w^A_m$ & $\,K\,$ &\mbox{ Form}&\mbox{ Degree}
\\ \hline
\endfirsthead
\caption[]{(continued)}\\ \hline $\phantom{n}m\phantom{n}$& $w_m$ & $\,K\,$ &\mbox{
Form}&\mbox{ Degree}
\\
\hline
\endhead
 \hline 3 & $12.8340776752\ldots$ &1& 3 & 10
 \\ \hline 33 & 5203280.6707141049\dots &2& 4 & 16
 \\ \hline 43 & 170133239.5931416562\dots &3& 3 & 18
 \\ \hline 52 & 3506589575.3297508814\dots &3& 4 & 20
 \\ \hline 53 & 4906979442.0645648056\dots &3& 4 & 20
 \\ \hline 54 & 6965642842.5492071202\dots &3& 4 & 20
 \\ \hline 62 & 95994610190.3413097554\dots &4& 3 & 22
 \\ \hline 63 & 131582414832.0133343262\dots &4& 3 & 22
 \\ \hline 64 & 182480513596.8192404599\dots &4& 3 & 22
 \\ \hline 65 & 257327059360.7694099942\dots &4& 3 & 22
 \\ \hline 72 & 2512477187944.7980147749\dots &4& 4 & 24
 \\ \hline 73 & 3382770986274.7717090200\dots &4& 4 & 24
 \\ \hline 74 & 4595841393803.8776672113\dots &4& 4 & 24
 \\ \hline 75 & 6325468915069.0951433926\dots &4& 4 & 24
 \\ \hline 76 & 8869434642969.6959845223\dots &4& 4 & 24
 \\ \hline 83 & 84893140132749.4433303610\dots &5& 3 & 26
 \\ \hline 84 & 113336757486228.6751081110\dots &5& 3 & 26
 \\ \hline 85 & 152780346952246.4470059715\dots &5& 3 & 26
 \\ \hline 86 & 208816888813525.1041282031\dots &5& 3 & 26
 \\ \hline 87 & 291111233760547.9215716098\dots &5& 3 & 26
 \\ \hline 88 & 417810824838712.4573733198\dots &5& 3 & 26
 \\ \hline 93 & 2101339201083448.8581013596\dots &5& 4 & 28
 \\ \hline 94 & 2765414429020562.7271354169\dots &5& 4 & 28
 \\ \hline 95 & 3664279182837636.4124486346\dots &5& 4 & 28
 \\ \hline 96 & 4903834748673690.9286701217\dots &5& 4 & 28
 \\ \hline 97 & 6656234865994295.8866755890\dots &5& 4 & 28
 \\ \hline 98 & 9219414667236866.7539438955\dots &5& 4 & 28
 \\ \hline 99 & 13154452413669110.8651309899\dots &5& 4 & 28
 \\ \hline 104 & 67128787164683544.4523776297\dots &6& 3 & 30
 \\ \hline 105 & 87735771631740197.9352556672\dots &6& 3 & 30
 \\ \hline 106 & 115452705604234426.2027477576\dots &6& 3 & 30
 \\ \hline 107 & 153434175553796665.9933019243\dots &6& 3 & 26
 \\ \hline 108 & 206786189887861292.3744524831\dots &6& 3 & 30
 \\ \hline 109 & 284301510798701394.0821487141\dots &6& 3 & 30
 \\ \hline 110 & 402436933410128596.4722278511\dots &6& 3 & 30
 \\ \hline 111 & 595836902535559024.3188262687\dots &6& 3 & 30
 \\ \hline 115 & 2101589801997124293.4926488050\dots &6& 4 & 32
 \\ \hline 116 & 2729686255505205518.5715664555\dots &6& 4 & 32
 \\ \hline 117 & 3568896545606734776.3884783889\dots &6& 4 & 32
 \\ \hline 118 & 4710751164134714024.1028237957\dots &6& 4 & 32
 \\ \hline 119 & 6302155682133310144.5379747181\dots &6& 4 & 32
 \\ \hline 120 & 8593437376247698936.3526858911\dots &6& 4 & 32
 \\ \hline 121 & 12046660361010217345.2333308572\dots &6& 4 & 32
 \\ \hline 122 & 17615670890375383613.3569976086\dots &6& 4 & 32
 \\ \hline 125 & 50364867085052138405.0507562145\dots &7& 3 & 34
 \\ \hline 126 & 64793425556670001358.5275796370\dots &7& 3 & 34
 \\ \hline 127 & 83682707580452106791.2097623968\dots &7& 3 & 34
 \\ \hline 128 & 108746230244493146881.2976319596\dots &7& 3 & 34
 \\ \hline 129 & 142584503530083915887.6303725657\dots &7& 3 & 34
 \\ \hline 130 & 189322485836170786553.6385557962\dots &7& 3 & 34
 \\ \hline 131 & 255885036015648050531.0181807358\dots &7& 3 & 34
 \\ \hline 132 & 354803789323254565081.7957039872\dots &7& 3 & 34
 \\ \hline 133 & 511236181565006961686.6864799904\dots &7& 3 & 34
 \\ \hline 134 & 784028296517640221309.3273031865\dots &7& 3 & 34
 \\ \hline 136 & 1541633134566897630439.9212931033\dots &7& 4 & 36
 \\ \hline 137 & 1973936094367255380618.5869006801\dots &7& 4 & 36
 \\ \hline 138 & 2536282160116347359248.0508958835\dots &7& 4 & 36
 \\ \hline 139 & 3277120814346441078193.1926393845\dots &7& 4 & 36
 \\ \hline 140 & 4269133305200474922172.0302255769\dots &7& 4 & 36
 \\ \hline 141 & 5625976046648898147250.0406685544\dots &7& 4 & 36
 \\ \hline 142 & 7535031448888363741834.1039332556\dots &7& 4 & 36
 \\ \hline 143 & 10327388458447120861753.9448151515\dots &7& 4 & 36
 \\ \hline 144 & 14645940302969237174125.6523056427\dots &7& 4 & 36
 \\ \hline 145 & 21920980271484567377147.0958895900\dots &7& 4 & 36
 \\ \hline 147 & 46712811159099187620845.3387817426\dots &8& 3 & 38
 \\ \hline 148 & 59565056612990122680415.8493932533\dots &8& 3 & 38
 \\ \hline 149 & 76177853174918430594445.5483347921\dots &8& 3 & 38
 \\ \hline 150 & 97906064392590213580625.1938449228\dots &8& 3 & 38
 \\ \hline 151 & 126755987190626475051892.3954775655\dots &8& 3 & 38
 \\ \hline 152 & 165815400461133042447385.7661426883\dots &8& 3 & 38
 \\ \hline 153 & 220073534234554737523908.4115945786\dots &8& 3 & 38
 \\ \hline 154 & 298116449057161906517212.6325878034\dots &8& 3 & 38
 \\ \hline 155 & 416023538135200368061705.0307527308\dots &8& 3 & 38
 \\ \hline 156 & 607718608181858421565493.8445181846\dots &8& 3 & 38
 \\ \hline 161 & \phantom{1}2904804593623115788824211.9075380435\dots\phantom{1} &8& 4 & 40
 \\
\hline
\end{longtable}

\section{\bf The proof for the case ${\bf m=43.}$}
Let $m=43,$ then $\alpha=20$ and the polynomials $R_k=R_k^{20,20}$ are the ultraspherical
polynomials orthogonal on $[-1,1]$ with weight $(1-t^2)^{20}$ and normalized by the
condition $R_k(1)=1$. In this section we evaluate the quantity
 \begin{equation}\label{v7_26}
w^A_{43}=w^A_{43}(1/2)=\inf\left\{\frac{f(1)+f(-1)}{f_0}: \
f\in{\cal F}^A_{43}\right\},\quad {\cal F}^A_{43}={\cal
F}^A_{43}(1/2) .
 \end{equation}

According to the table of the results, the extremal function is an even eighteenth-degree
polynomial of the form 3 with $K=3;$ its fourteenth and sixteenth coefficients in
representation (\ref{v3}) with respect to the ultraspherical polynomials are zero. Thus,
the extremal polynomial has the form
 \begin{equation}\label{f^*}
 f^*(t)=\left(t^2-\frac{1}{4} \right)\cdot\prod\limits_{i=1}^3(t^2-a_i^2)^2\cdot(t^4+q\,t^2+r), \quad
 f^*_{14}=f^*_{16}=0.
 \end{equation}
We are to choose parameters of the polynomial, i.e., the zeros $\{a_i\}$ and the
coefficients $q,\, r;$ in addition, the following two conditions must be satisfied:

(c1) the (double) zeros $a_i,\ 1\le i\le 3,$ lie on $\left(0,1/2\right);$

(c2) the polynomial $t^4+q\,t^2+r$ is nonnegative on $[0,1/2]$ (and, in reality, it will
be positive on all the axis).

Simultaneously with the function (\ref{f^*}) we are to construct the quadrature formula
 \begin{equation}\label{F26}
 f_0=\frac{1}{\vartheta(\alpha)}\int\limits_{-1}^{1}f(t)(1-t^2)^\frac{m-3}{2}dt=
 L(f)-\gamma_{14} f_{14}-\gamma_{16} f_{16},
 \end{equation}
\begin{equation}\label{L26}
L(f)= \sum\limits_{i=1}^{3} A_i\, f(a_i) + A_0\,
f\left(\frac{1}{2}\right)+ A_4\, f\left(1\right),
\end{equation}
with the following properties:

(p1) the formula is sharp on the set ${\cal P}^A_{18}$ of even algebraic polynomials of
degree up to 18;

(p2) $L(R_{2k})\ge 0, \ k\ge 0;$ in addition (as a consequence of the preceding
condition),   $L(R_{2k})= 0,\  1\le k\le 9,\ k\ne 7,\,8;$

 (p3) the coefficients $A_i,\ 0\le i\le 4,$ of  formula
(\ref{L26}) are nonnegative.

The number of all the parameters here is twelve. Namely, there are five coefficients
$A_i,$ three unknowns $a_i,$ two coefficients of the polynomial $t^4+q\,t^2+r,$ and two
coefficients $\gamma_{14},$ $\gamma_{16}.$ In order to find these parameters, we
construct the set ${\Sigma}_{12}$ of twelve (nonlinear) equations. Ten equations are
given by the condition that the quadrature formula (\ref{F26}) is sharp for even
algebraic polynomials of the eighteenth degree, i.e., the condition (p 1); and two
equations are given by the conditions that the coefficients $f^*_{14},\ f^*_{16}$ of
representation of the unknown polynomial  $f^*$ by the system $\{R_{2k}\}$ are equal to
zero.

In order to simplify ${\Sigma}_{12}$, we use the following considerations. We use a basis
of polynomials in the space ${\cal P}^A_{18}$ such that after substitution of these
polynomials in the quadrature formula (\ref{F26}) we get the simplest equations. In
addition, we replace the squared zeros $a_i,\  1\le i\le 3,$ of the polynomial
(\ref{f^*}) by their symmetric functions $U_i,\ 1\le i\le 3.$ To be exact, we find the
polynomial \mbox{$P_3(t)=\prod\limits_{i=1}^3(t^2-a_i^2),$} whose square is in the
expansion (\ref{f^*}), in the form
$$ P_3(t)= t^6-U_{2}t^{4}+U_{1}t^{2}-U_{0}.$$ The
system of equations $\Sigma_{12}$ constructed in this way has several solutions, and only
one of them satisfies  conditions  (c1), (c2), (p2), (p3).

One can find the construction of the system of equations and its solution in the proof of
the following theorem and Lemma~\ref{lb2}. We solved the system with the help of the
Maple package of analytic computations. We do not give these constructions and
computations. We give a result of  the computations, i.e., we give the concrete function
(\ref{f^*}) and the quadrature formula (\ref{F26}) + (\ref{L26})  and prove that they
solve the problem.

To state the results of this part of the work, we need some definitions and notation. Let
$H$ be the following polynomial of the third degree:
\begin{eqnarray}\label{f2}
H(z) = {z}^{3}-{\frac {1835489}{2079100}}\,{z}^{2} + {\frac
{590059779}{1287046064}}\,z-{\frac {106321508304907}{
2129617205027920}}.
\end{eqnarray}
The polynomial has one real and two complex zeros:
\begin{eqnarray*}
 z_1 & =&     0.1411134854416294\ldots,\\
 z_{2,3} & =& 0.3708575711650012\ldots\pm i \cdot 0.4650367196476447\ldots\,.
\end{eqnarray*}

 Let us agree to denote the zero  $z_1$ of the polynomial $H$ by $\xi,$ so that
 \begin{equation}\label{f3}
 \xi=z_1=0.1411134854416294\ldots\,.
 \end{equation}
Let
 \begin{eqnarray}
 \nonumber q& = & 2\,\xi - {\frac {179}{100}} = -1.5077730291167411\ldots,
 \\
 \nonumber r&=& 3\,{\xi}^{2} -{\frac
{3914589}{1039550}}\,\xi+{\frac {40170654239}{32176151600}} =
0.7768145246622512\ldots,
 \end{eqnarray}
 \begin{eqnarray} \label{zetadef}
 \zeta_0 & = & {\frac {779}{1018759}}\,\xi-{\frac {5570}{53994227}}\, , \nonumber \\
 \zeta_1 & = & {\frac {1930}{20791}}\,\xi-{\frac{10605}{1101923}}\, , \\
 \zeta_2 & = & \xi\, . \nonumber
 \end{eqnarray}

The numbers $q, r,\zeta_0,\zeta_1,\zeta_2$ define the polynomial of  the eighteenth
degree
 \begin{equation}\label{f7}
 f^*(t)=\left(t^2-\frac{1}{4}\right)\left(t^6-\zeta_2\,t^4+\zeta_1\,t^2-\zeta_0 \right)^2
 (t^4+q\,t^2+r).
 \end{equation}
We denote by $f^*_k$ its Fourier coefficients  in the expansion with respect to the
polynomials $\{R_k\}$. The polynomial $f^*$ is an even polynomial. Therefore its
coefficients with odd indexes are equal to zero. Hence the expansion has the following
form
 \begin{equation}\label{f8}
 f^*(t)=\sum_{k=0}^{18} f_{2k}^* R_{2k}(t)
 \end{equation}
Finally, let
 $\{ a_i \}_{i=1}^{3}$ be  positive zeros of the polynomial
 \begin{equation}\label{f9}
 g(z)=z^6-\zeta_2\,z^4+\zeta_1\,z^2-\zeta_0;
 \end{equation}
they all belong to $(0,1/2)$ and have the following approximate values:
 \begin{eqnarray}\label{f10}
 a_1 & = & 0.0380726602850886\ldots\,, \nonumber \\
 a_2 & = & 0.1725867591939439\ldots\,, \\
 a_3 & = & 0.3314781569445825\ldots\,. \nonumber
 \end{eqnarray}
By (\ref{f7}), $\pm a_i,\ 1\le i\le 3,$ are the double zeros of the function $f^*.$

The  following statement is the main result of this section.

 \begin{Thm}\label{tb1}
The polynomial $f^*,$ defined by $(\ref{f2})-(\ref{f7}),$ belongs to the set ${\cal
F}_{43}^A$ and is the unique $($to within a positive constant factor$)$ extremal function
$($solution$)$ of the problem $(\ref{v7_26}).$ Moreover, we have
 \begin{equation}\label{f11}
 w^A_{43}=\frac{2 f^*(1)}{f_0^*}= 170133239.5931416562399728\ldots\ .
 \end{equation}
 \end{Thm}

The proof of Theorem~\ref{tb1} requires some preliminary statements.

For the ultraspherical polynomials $R_n=R_n^{\alpha,\alpha},\ \alpha=(m-3)/2,$ for any
$m\ge 2,$ in particular, for $m=43,$  the following recurrence relation holds (see, for
instance, \cite[p.64, formula (4.5)]{Lev-83}):
  $$
  R_{n+1}(t) =  \frac{(2n+m-2)\, t\, R_{n}(t)-n\, R_{n-1}(t)}{n+m-2}, \quad n\geq1,
  $$
$$
   R_0(t) = 1,\ \ R_1(t)=t.
$$
This relation is usable for evaluation of the coefficients of the polynomials $R_n,\ n\ge
1,$ in their expansion by degrees of $t.$ In the sequel on several occasions we represent
a polynomial \mbox{$f(t)=\sum_{k=0}^{\nu} c_k(f)t^k$} of degree $\nu$ by the system
$\{R_k\}:$
\begin{equation}\label{fR}
f(t)=\sum_{k=0}^{\nu}f_kR_k(t).
\end{equation}
We do this by the following known scheme. Suppose $c_{\nu}(R_{\nu})$ is the leading
coefficient of the polynomial $R_{\nu}.$ Then for the expansion (\ref{fR}) we have
$f_{\nu}=c_{\nu}(f)/c_{\nu}(R_{\nu}).$ The degree of the polynomial
$f^1=f-f_{\nu}\,R_{\nu}$ is equal to $\nu-1.$ Repeating this process for the polynomial
$f^1,$ we get the coefficient $f_{\nu-1}$ and so on.

With the help of the polynomial (\ref{f9}) we define the polynomials
\begin{equation}\label{la_defg}
 g^i(t)=(t^2-1/4)(t^2-1)\, g(t) / (t^2-a_i^2) \quad i=1,\ 2,\  3.
\end{equation}
In what follows we use the notation
\begin{equation}
 g^i(t)=\sum_{j=0}^{9}g^{i}_{2j}\,R_{2j}(t)
\end{equation}
for coefficients in representations of these polynomials by the system of ultraspherical
polynomials. Let us introduce the quantities
\begin{equation}\label{la_def}
\lambda(a_i)=g^i_0/g^i(a_i), \quad i=1,\ 2,\ 3.
\end{equation}
Their numerical values are the following:
\begin{eqnarray}\label{h35}
 \lambda(a_1) & = & 0.4835866972149467\ldots\,,\nonumber\\
 \lambda(a_2) & = & 0.4319241073046564\ldots \,,\\
 \lambda(a_3) & = & 0.0815919008616610\ldots \,. \nonumber
\end{eqnarray}
Let also
\begin{equation}\label{h33}
\lambda(1)=-{\frac {1}{23009085}}\,\left({\frac
{214703\,\xi-24075}{24221\,\xi-26423}}\right)
=0.00000001175549\ldots,
\end{equation}
\begin{equation}\label{h34}
\lambda\left(\frac{1}{2}\right)={\frac
{928514048}{90945}}\,\left({\frac {53\,\xi-15}{138423916\,
\xi-46036387}}\right)
 = 0.002897282863\ldots.
\end{equation}

The following statement contains the quadrature formula on the set of functions
(\ref{h0A}) for $m=43.$ The proof of Theorem~\ref{tb1} is based on this formula.

 \begin{Lemma}\label{lb2}
For functions $f(t)=\displaystyle\sum_{j=0}^\infty f_{2j} R_{2j}(t) \in \Phi^A_{43}$ we
have the following quadrature formula:
\begin{equation}\label{h31}
 f_0={\frac {1412926920405}{549755813888}}
\int _{-1}^1 f(t)(1-t^2)^{\frac{43-3}{2}}dt = L(f)-\sum_{\nu\ge 1}
 L(R_{2\nu})f_{2\nu};
 \end{equation}
here $L$ is the functional
 \begin{equation}\label{h32}
 L(f)=\lambda(1)f(1)+\lambda\left(\frac{1}{2}\right)f\left(\frac{1}{2}\right)
 +\sum_{i=1}^{3}\lambda(a_i)f(a_i)
  \end{equation}
with  coefficients defined by {\rm(\ref{la_def})}---{\rm(\ref{h34})}. The functional $L$
has the following properties:
 \begin{equation}\label{h390}
 L(1)=1;
 \end{equation}
 \begin{equation}\label{h38}
 L(R_{2\nu})=0,\quad  \nu= 1,\ 2,\ldots,\ 5,\ 6,\ 9;
 \end{equation}
 \begin{equation}\label{h39}
 L(R_{2\nu})>0,\quad  \nu= 7,\ 8,\ \nu\ge 10.
 \end{equation}
 \end{Lemma}

 \noindent {P\ r\ o\ o\ f.\ }
 In the proof of the lemma we use the ideas similar to those
applied in the proof of Lemma~4.2 from \cite{Arestov and Babenko} in the investigations
of the quantity $w_4(1/2).$

We first find  what conditions on the real nodes
 \begin{equation}\label{h40}
 0<A_1<A_2<A_3<\frac{1}{2}
 \end{equation}
and the coefficients $\lambda_1, \ \lambda_{1/2}, \ \lambda_{A_1}, \ \lambda_{A_2},\
\lambda_{A_3}, \ \gamma_{14}, \ \gamma_{16}$ must hold  if the quadrature formula
 \begin{equation}\label{h41}
 f_0={\frac {1412926920405}{549755813888}}
\int_{-1}^1 f(t)\cdot(1-t^2)^{\frac{43-3}{2}}
  dt = {\cal L}(f)-\gamma_{14}\,f_{14}-\gamma_{16}\,
 f_{16},
 \end{equation}
 \begin{equation}\label{h42}
 {\cal L}(f)=\lambda_1 f(1)+\lambda_{1/2}f\left(\frac{1}{2}\right) +\sum_{i=1}^{3}
 \lambda_{A_i} f(A_i),
 \end{equation}
is valid on the set of all even polynomials $f(t)=\sum_{k=0}^{9} f_{2k} R_{2k}(t)$ of
degree up to 18. Next we add several conditions that together with the previous ones
become sufficient for construction of the required quadrature formula on the set
$\Phi^A_{43}.$

Let
 \begin{equation}\label{h45}
 \chi(t)= \prod^{3}_{i=1} (t^2-A_i^2) = t^6-U_2\,t^4+U_1\,t^2-U_0.
 \end{equation}
We consider  the following polynomial of the tenth degree:
 \begin{equation}\label{h46}
 \sigma(t)=(t^2-1)\left(t^2-\frac{1}{4}\right)\chi(t).
 \end{equation}
With the help of the polynomial, we define some more algebraic polynomials  (of degree 
up to 18). Let  $\varphi^1=\sigma$ and
 $$
 \varphi^1(t)=\sum_{k=0}^{9}
 \varphi^1_{2k} R_{2k}(t)
 $$
be the expansion of the polynomial by the system $\{R_k\}^{\infty}_{k=0}.$ In fact, the
degree of the polynomial is  10 and therefore, $\varphi^1_{14}=\varphi^1_{16}=0.$ Our
interest is the coefficient of $R_0$ in the expansion. For this coefficient we have
$$
 \varphi^1_0=-{\frac {29}{141470}}\,U_{{2}}+{\frac {49}{12126}}\,U_{{1}}-{\frac {
287}{1290}}\,U_{{0}}+{\frac {23}{1442994}} .
 $$
Similarly, for the polynomial
 $$\varphi^2(t)=t^2\sigma$$
we have
 \begin{eqnarray*}
 \varphi^2_0 & = & -{\frac {23}{1442994}}\,U_{{2}}+{\frac {29}{141470}}\,U_{{1}}-{\frac {
49}{12126}}\,U_{{0}}+{\frac {1}{642678}} .
 \end{eqnarray*}
Substituting the polynomials $\varphi^1,\,\varphi^2$ in (\ref{h41}), we get the first
necessary condition of existence of the formula:
 $$\varphi^1_0=0,\ \  \varphi^2_0=0.$$
The last relations give the following system of two linear equations with three unknowns:
 {
 \begin{eqnarray*}
 &\phantom{+}&-{\frac {29}{141470}}\,U_{{2}}+{\frac {49}{12126}}\,U_{{1}}-{\frac {
287}{1290}}\,U_{{0}}+{\frac {23}{1442994}}=0,
 \\&\phantom{+}&-{\frac {23}{1442994}}\,U_{{2}}+{\frac {29}{141470}}\,U_{{1}}-{\frac {
49}{12126}}\,U_{{0}}+{\frac {1}{642678}}=0;
 \end{eqnarray*}
}this is equivalent to
 \begin{eqnarray} \label{h47}
 U_0 & = & {\frac {779}{1018759}}\,U_{{2}}-{\frac {5570}{53994227}}, \nonumber \\
 U_1 & = & {\frac {1930}{20791}}\,U_{{2}}-{\frac {10605}{1101923}}.
 \end{eqnarray}

We consider the polynomial
 $\varphi^3(t)=t^4\cdot\sigma(t);$  it is of the fourteenth degree. For
the coefficients in the expansion  of the polynomial with respect to the system of the
ultraspherical polynomials, we have
 $$ {\varphi^3_{0}} = -{\frac {1}{642678}}\,U_{{2}}+{\frac
{23}{1442994 }}\,U_{{1}}-{\frac {29}{141470}}\,U_{{0}}+{\frac
{3}{18209210}},
 $$
 $$\varphi^3_{14}={\frac {847872}{4581527}} ,\ \ \varphi^3_{16}=0.$$
Substituting the polynomial $\varphi^3$ in (\ref{h41}), we get
\begin{equation}
\varphi^3_{0}=-\gamma_{14}\,\varphi^3_{14};
\end{equation}
this is equivalent to
\begin{eqnarray}\label{gamma12}
 \gamma_{14}&=& {\frac {4581527}{544908681216}}\,U_{
{2}}-{\frac {4581527}{53194530816}}\,U_{{1}}+ \\ &+&{\frac
{132864283}{ 119948451840}}\,U_{{0}}-{\frac
{4581527}{5146359767040}}.\nonumber
\end{eqnarray}
Similarly, using the polynomial of the sixteenth degree $\varphi^4(t)=t^6\cdot\sigma(t)$
and the expression (\ref{gamma12}) for $\gamma_{14}$, we find the coefficient
$\gamma_{16}:$
\begin{eqnarray}\label{gamma13}
 \nonumber\gamma_{16} &=& {\frac
{325288417}{24323459973120}}\,{U_{{2}}}^{2}-{\frac {7481633591}{
54613051637760}}\,U_{{1}}U_{{2}}+\\
 &+&{\frac
{9433364093}{5354220748800}}\, U_{{0}}U_{{2}}-{\frac
{543272890133}{11577966947205120}}\,U_{{1}}-\\ &-&{\frac
{19925060923 }{31207458078720}}\,U_{{0}}-{\frac
{114538175}{19458767978496}}\,U_{{2 }}+{\frac
{39304679}{78761679912960}}.\nonumber
\end{eqnarray}

Let us introduce the polynomials
 \begin{eqnarray*}
 h_{i}(t) & = & \sigma(t)/(t^2-(a_i)^2), \quad 1\leq i \leq3,\\
 h_{4}(t) & = & \left(t^2-\frac{1}{4}\right)\chi(t)=\sigma(t)/(t^2-1),\\
 h_{5}(t) & = & (t^2-1)\chi(t)=\sigma(t)/(t^2-1/4).
 \end{eqnarray*}
Substituting the polynomials $h_4,\, h_5$ in (\ref{h41}), we get
 \begin{eqnarray*}
 \lambda_1 & = & -{\frac {1}{636615}}\,{\frac {189\,U_{{2}}-3619\,U_{{1}}+192465\,U_{{0
}}-15}{U_{{2}}-U_{{1}}+U_{{0}}-1}},
\\
 \lambda_{1/2} & = & {\frac {512}{636615}}\,{\frac {147\,U_{{2}}-2303\,U_{{1}}+103635\,U_{{0
}}-15}{4\,U_{{2}}-16\,U_{{1}}+64\,U_{{0}}-1}}.
 \end{eqnarray*}
Substituting $U_0,\ U_1$ from (\ref{h47}) in formulas for $\lambda_1$, $\lambda_{1/2},$
we get
 \begin{eqnarray}\label{lambda1}
 \lambda_1 & = & -{\frac {1}{23009085}}\,{\frac {214703\,U_{{2}}-24075}{24221\,u
_{{2}}-26423}},\\
 \lambda_{1/2} & = & {\frac {928514048}{90945}}\,{\frac {53\,U_{{2}}-15}{
138423916\,U_{{2}}-46036387}} .
 \end{eqnarray}

In the same way, with the help of the polynomials $h_i,\ 1\le i\le 3,$ we find
 \begin{equation}\label{lambda_ai}
 \lambda_{A_i}=h_0^i/h^i(A_i),\  i=1,\ 2,\ 3.
 \end{equation}
We do not give explicit formulas  for the coefficients $\lambda_{A_i},$ because the
expressions  occupy too much room.

Substituting the polynomials
 $\varphi^5(t)=t^8\sigma(t)$ in (\ref{h41}) and considering  the
formulas obtained above  for $ \gamma_{14}, \gamma_{16}, U_0 , U_1,$ we get the following
equation:
 \begin{equation}\label{syst}
  F(U_2)=0,
 \end{equation}
where
 \begin{equation}\label{ur1}
 F(U_2)={U_{{2 }}}^{3}-{\frac {1835489}{2079100}}\,{U_{{2}}}^{2}+{\frac {590059779}{
1287046064}}\,U_{{2}}-{\frac {106321508304907}{2129617205027920}}.
 \end{equation}

The polynomial $F(U_2)$ coincides with the polynomial $H$. Thus, the conditions
$H(U_2)=0$ and (\ref{h47}), and (\ref{gamma12})---(\ref{lambda_ai}) are the necessary
conditions for the existence of the quadrature formula~(\ref{h41}).

As mentioned above, the polynomial $H$ has one real and two complex zeros; in (\ref{f3})
we denoted its first real zero  by $\xi$. From now on we suppose that
 $$
 U_2=\zeta_2=\xi,\ U_1=\zeta_1,\ U_0=\zeta_0,
 $$
where $\zeta_0,$ $\zeta_1,$ $\zeta_2$ are defined by (\ref{zetadef}). Under these
assumptions,  the polynomial $\chi$ defined by (\ref{h45}) coincides with the polynomial
(\ref{f9}).  Consequently,
 \begin{equation}\label{h69}
 A_i=a_i, \quad i=1,\ 2,\ 3.
 \end{equation}
This implies that the coefficients
 $\lambda_1,\lambda_{1/2},\lambda_{A_1},\lambda_{A_2},\lambda_{A_3}$ of the functional $\cal L$  (see (\ref{h42})) coincide with the
coefficients (\ref{la_def}) --- (\ref{h34}) of the functional $L$ defined by (\ref{h32}).
Therefore, (\ref{h41}) takes the form
 \begin{equation}\label{h70}
 f_0={\frac {1412926920405}{549755813888}}\int_{-1}^1 f(t)(1-t^2)^\frac{43-3}{2}dt
 = L(f)-\gamma_{14} f_{14}-\gamma_{16} f_{16}.
 \end{equation}

Now we can assert that the formula is valid for the polynomials $h^1,\ h^2,\ \ldots,\
h^5,$ $\varphi^1,\ \varphi^2,\ \ldots,\ \varphi^5.$ These ten polynomials form the basis
in the set ${\cal  P}^A_{18}$ of even polynomials of degree up to eighteen. Therefore,
the quadrature formula (\ref{h70}) holds for any polynomial from ${\cal P}^A_{18}.$
Substituting the polynomials $R_{14}, \ R_{16}$  in the formula, we obtain
 \begin{equation}\label{h55}
 \gamma_{14}=\,L(R_{14}),\quad \gamma_{16}=\,L(R_{16}).
 \end{equation}
Hence, (\ref{h70}) can be written in the form
 \begin{equation}\label{h71'}
 f_0={\frac {1412926920405}{549755813888}}\int_{-1}^1 f(t)(1-t^2)^{\frac{43-3}{2}}dt = L(f)-\sum_{\nu=7}^{8} L(R_{2\nu})
 f_{2\nu},
 \end{equation}
 where
 \begin{equation}\label{h71}
 L(R_\nu) = \lambda(1)+
 \lambda\left(\frac{1}{2}\right) R_\nu\left(\frac{1}{2}\right)+
 \sum_{i=1}^{3}\lambda(a_i)R_\nu(a_i),
 \end{equation}
 $a_1,\ a_2,\ a_3$ are positive zeros of the polynomial given by (\ref{f9}).

Hence, in particular, $L(R_{2\nu})=0$ for $\nu= 1,\, 2,\ldots,\, 5,\, 6,\, 9,$  $\nu \neq
7,\, 8,$ i.e., the property (\ref{h38}) holds. Substituting the polynomial  $f(t)\equiv
1$ in (\ref{h71'}), we also get (\ref{h390}).

Obviously, (\ref{h71'}) can be extended to all the class of functions $\Phi^A_{43},$ if
we write this formula in the form~(\ref{h31}).

In order to complete the proof, it remains to verify inequalities (\ref{h39}). To
validate these inequalities, we use the following two properties of the ultraspherical
polynomials $\{R_k\}$:

 (1)~$R_\nu(1)=1,\ R_\nu(-1)=(-1)^\nu,\  \nu\ge 0, \ $

 (2) on $(-1,1)$, the polynomials $R_\nu$ converge pointwise to zero as
$\nu\to\infty,$ to be exact, the estimates (\ref{ef_f}) hold for the polynomials.

For any $\nu\ge 0$ we have {\small
 \begin{eqnarray*}
 L\left(R_\nu\right) & \ge & \lambda(1)+
 \left(\lambda\left(\frac{1}{2}\right)+\sum_{i=1}^{3}\lambda(a_i)\right)
 \min\left\{R_\nu\left(\frac{1}{2}\right),R_\nu(a_1),\ldots,R_\nu(a_3)\right\}=
 \\ & = & \lambda(1) + (1-\lambda(1))
 \min\left\{R_\nu\left(\frac{1}{2}\right),R_\nu(a_1),R_\nu(a_2),R_\nu(a_3)\right\}
 \geq
 \\ & \geq & \lambda(1)- (1-\lambda(1))
 \max\left\{\left| R_\nu\left(\frac{1}{2}\right)\right| ,
 \left|R_\nu(a_1)\right|,\left|R_\nu(a_2)\right|,\left|R_\nu(a_3)\right|\right\}.
 \end{eqnarray*}
 }
Let us use the following efficient estimates for the ultraspherical polynomials
\begin{equation}\label{ef_f}
\left| R_n^{\alpha,\alpha}(t)\right| \leq \frac{{A}(n,m)}{(1-t^2)^\gamma},\quad
\alpha=(m-3)/2,\quad -1<t<1,
 \end{equation}
$$ {n\geq \max\{3,m-4\},}\ \  m\geq 4, $$
\begin{equation}\label{ef_Anm}
{A}(n,m)=\ {\Gamma\left(\frac{m-1}{2}\right)} \frac{
\sqrt{2}(2+\sqrt{2})^{m-4}}{(n+1)^{\frac{m-2}{2}}},\quad \gamma=\frac{m-2}{4}.
\end{equation}
Estimates of the Jacobi polynomials and, particularly, of the ultraspherical polynomials
has a rich history; the estimate (\ref{ef_f}) is contained in (\cite{Shtrom1},
Lemma~2.1).

By the assumption of the lemma, the estimate (\ref{ef_f}) can be used for  $n\ge 39.$ Let
us recall that $ a_i \in (0,1/2).$ Therefore, for any $\nu\ge 1$ we have
 \begin{eqnarray}
\nonumber & \max & \left\{\left|
R_\nu\left(\frac{1}{2}\right)\right| ,
 \left|R_\nu(a_1)\right|,\left|R_\nu(a_2)\right|,\left|R_\nu(a_3)\right|\right\} \leq
 \\ & \leq &
 \frac{\Gamma(21)\sqrt{2}(2+\sqrt{2})^{39}\,}
 {(1-1/4)^{\frac{41}{4}}\cdot(\nu+1)^{\frac{41}{2}}}\leq
  \frac{10^{42}}{(\nu+1)^{\frac{41}{2}}}.\nonumber
 \end{eqnarray}
It follows easily that  $L\left(R_\nu\right)>0$ for
 $\nu\ge 2000$. For other $\nu$ the conditions (\ref{h39}) can be verified by
direct computation with the use of the  Maple program. This completes the proof of
Lemma~\ref{lb2}.

 \begin{Lemma}\label{lb1}
The function $f^*$ defined by $(\ref{f2})$--$(\ref{f7})$ belongs to the set~${\cal
F}^A_{43}.$
 \end{Lemma}

 \noindent {P\ r\ o\ o\ f.\ }
The polynomial $\pi_2(t)=t^4+q\,t^2+r$ in the right-hand side of (\ref{f7}) is positive
on $[-1,1]$. Therefore,
 \begin{equation}\label{h25}
 f^*(t)\le 0,\quad  -\frac{1}{2}\le t\le\frac{1}{2}.
 \end{equation}

It remains to prove that the coefficients  $f_{2k}^*$ in the representation (\ref{f8}) of
the function (\ref{f7}) are nonnegative; moreover, $f^*_0>0.$ We evaluate the
coefficients of the representation
 $$
 f^*(t)=\left(t^2-\frac{1}{4}\right)\left(t^6-\zeta_2\,t^4+\zeta_1\,t^2-\zeta_0 \right)^2(t^4+q\,t^2+r)
$$ with respect to the system of ultraspherical polynomials:
 \begin{eqnarray}\label{h29}
 f_{18}^* & = &
 {\frac {439025664}{6248961695}},
 \\ f_{16}^* & = & f_{14}^*=0\,,\nonumber
 \\ f_{12}^* & = & 0.0417636211579982720349265\ldots\,,\nonumber
 \\ f_{10}^* & = & 0.0292289140950023331221331\ldots\,,\nonumber
 \\ f_{8}^* & = &  0.0078138941286457485943531\ldots\,,\nonumber
 \\ f_{6}^* & = &  0.0009446460463425510583249\ldots\,,\nonumber
 \\ f_{4}^* & = &  0.0000488080874818309645154\ldots\,,\nonumber
 \\ f_{2}^* & = &  0.0000008242017491816358353\ldots\,,\nonumber
 \\ f_{0}^* & = &  0.0000000017639878907626135\ldots\,,\nonumber
 \end{eqnarray}
Thus, the representation (\ref{v3}) of the function $f^*$ has nonnegative coefficients
and \mbox{$f_0^*>0.$} This completes the proof of Lemma~\ref{lb1}.

 \ \

The following lemma will allow us to show that the function  $f^*$ is the unique (to
within a positive constant factor) extremal function of the problem~(\ref{v7_26}).

 \begin{Lemma}\label{lb3}
Suppose that $f$ is an even polynomial of degree up to eighteen with the following
properties:

 $(1)$ the representation
 \begin{equation}\label{E1}
 f(t)=h(t)e(t)
 \end{equation}
holds, where
 $$
 h(t)= \left(t^2-\frac{1}{4}\right)\prod_{i=1}^{3}(t^2-a_i^2)^2,
 $$
 and $e$ is some even polynomial of degree up to four,

 $(2)$ in the representation
 \begin{equation}\label{E2}
 f(t)=\sum_{{\nu}=0}^{18}f_{\nu} R_{\nu}(t)
 \end{equation}
of the polynomial $f$ by the polynomials $R_k,$ the fourteenth and sixteenth coefficients
are equal to zero:
 \begin{equation}\label{E3}
 f_{14}=0,\quad f_{16}=0.
 \end{equation}
 \noindent Then the polynomial $f$ coincides to within a constant factor with the
function {\rm(}polynomial{\rm)} $f^*,$ i.e., $f=cf^*,$ where $c={\rm const};$ in
addition, if $f\in {\cal F}_{43}^A,$  then $c>0.$
 \end{Lemma}

 \noindent {P\ r\ o\ o\ f.\ }
 The proof of the lemma will be divided into several steps.

(1) We first suppose that the polynomial $e$ in the representation (\ref{E1}) has the
form
 \begin{equation}\label{hg}
 e(t)=t^4+e_1t^2+e_0.
 \end{equation}
Let us write out the explicit expressions for the coefficients $f_{14}$ and $f_{16}$ in
the representation (\ref{E2}) of the polynomial~$f.$ It is easily shown that
 \begin{equation}\label{E4}
 f(t)=f_{18}R_{18}(t)+f_{16}R_{16}(t)+f_{14}R_{14}(t)+\varphi(t),
 \end{equation}
where
 $$
 f_{18}={\frac {439025664}{6248961695}},
 $$
 $$
 f_{16}=-{\frac {75694080}{325288417}}\,\xi+{\frac {37847040}{325288417}}\,e_{
{1}}+{\frac {338731008}{1626442085}} ,
 $$
 \begin{eqnarray*}
 f_{14} & = & -{\frac {3372580030464}{6763071477847}}\,\xi+{\frac {847872}{4581527}}
\,{\xi}^{2}+{\frac {847872}{4581527}}\,e_{{0}} + \\& + &{\frac
{86694912}{ 325288417}}\,e_{{1}}-{\frac
{1695744}{4581527}}\,e_{{1}}\xi+{\frac {
6437476926996480}{26166323547790043}},
 \end{eqnarray*}
and, finally, $\varphi$ is an even polynomial of the twelfth degree.

Let the function $f$ satisfy the conditions (\ref{E3}), i.e., let the right-hand sides of
the last two relations be equal to zero. As a result, we obtain a system of two linear
equations with two unknowns $e_1 $ and $ e_0.$ The system determinant $\Delta$ has the
form 
$$
\Delta=-{\frac {52381392652633374720}{1299065844351175880963}}+{\frac {
3916313168896327680}{24510676308512752471}}\,\xi-{\frac {2595980574720
}{16604338082711}}\,{\xi}^{2}.
$$
It follows easily that by (\ref{f3}), $\Delta=-0.0208885275\ldots\ne 0.$ Consequently,
the system under consideration, i.e., the system of two conditions (\ref{E3}), has the
unique solution $e_0,\, e_1.$ By (\ref{f7}), the function $f^*$ has the representation
(\ref{E1}), (\ref{hg}) and by (\ref{h29}), the conditions (\ref{E3}) hold for this
function. Hence,  $e_1=q$ and $e_0=r.$ Thus, $f=f^*$ is a unique function satisfying the
requirements (\ref{E1}), (\ref{hg}) and (\ref{E3}).

(2) Suppose now that the polynomial $e$ in the representation (\ref{E1}) of the function
$f$ has the form $e(t)=ct^4+e_1\,t^2+e_0$ and $c\ne 0.$ Then the function
$f_c=\ds\frac{f}{c}$ satisfies all the conditions of the lemma and the second factor in
the representation (\ref{E1}) of the function has the form (\ref{hg}). It follows, as we
have shown just now, that $\ds\frac{f}{c}=f^*$ or, equally, $f=cf^*.$ Clearly, if $f\in
{\cal F}_{43}^A,$ then $c>0.$

 (3) Finally, let us show that if the function $f$ satisfies the assumptions of the lemma
and the second factor in (\ref{E1}) has the form $e(t)=e_1\,t^2+e_0,$  then $f\equiv 0.$
Indeed,  the function $\overline f=f^*+f$ satisfies all the assumptions of the lemma. For
this function the formula $\overline f=h\overline e$ is valid, where $\overline e$ is the
polynomial of the fourth degree  $\overline e(t)=t^4+\overline e_1\,t^2 +\overline e_0$
with the coefficients $ \overline e_1  =q+e_1, \ \overline e_0=r+e_0.$ The polynomial
$\overline e$ has the form (\ref{hg}). Hence, as we have already proved, $\overline
f=f^*$ and therefore $f\equiv 0.$ This completes the proof of Lemma~\ref{lb3}.

     \ \

 {P\ r\ o\ o\ f\ \ of Theorem~\ref{tb1}.} For the function $f^*$ defined by  $(\ref{f7}),
(\ref{f8})$ we have
 $$
 f^*\left(\frac{1}{2}\right)=f^*(a_1)=f^*(a_2)=f^*(a_3)=0, \quad
 f^*_{14}=f^*_{16}=0.
 $$
Therefore, for this function, inequality (\ref{h31}) turns into equality
 $$
 f^*_0=\lambda(1)f^*(1).
 $$
 This implies that
 $$
 \frac{f^*(1)}{f^*_0}=\frac{1}{\lambda(1)}=-23009085\,\left({\frac {24221\,\xi-26423}{214703\,\xi-24075}}\right).
 $$
 By Lemma~\ref{lb1}, the function $f^*$ belongs to the set ${\cal F}_{43}^A$.  Therefore,
  $$
 w^A_{43}\le\frac{2}{\lambda(1)}=\frac{2\cdot f^*(1)}{f^*_0}.
 $$
Thus, in order to prove (\ref{f11}), it is sufficient to show that
 \begin{equation}\label{h72}
 w^A_{43}\ge\frac{2}{\lambda(1)}.
 \end{equation}

For any function $f\in {\cal F}^A_{43}$ the quadrature formula
 (\ref{h31})---(\ref{h32}) is valid. The coefficients
$ \lambda(1), \lambda\left(\ds\frac{1}{2}\right), \lambda(a_1), \lambda(a_2),
\lambda(a_3)$ of the formula are positive, and the coefficients $L(R_{2\nu}), \ \nu\ge
1,$ are nonnegative. By  properties (\ref{v3})---(\ref{v6}) of the function $f\in {\cal
F}^A_{43}$, we have the following estimate:
 \begin{equation}\label{h73}
 f_0=L(f)-\sum_{\nu\ge 1} L(R_{2\nu})f_{2\nu} \le \lambda(1)f(1).
 \end{equation}
 This implies inequality (\ref{h72}). This proves the assertion
(\ref{f11}) and simultaneously  proves that the function $f^*$ is extremal in the problem
(\ref{v7_26}).

In order to prove the theorem, it remains to check  that the function  $f^*$ is the
unique, to within a positive constant factor, extremal function. The even function
 \begin{equation}\label{f-extr}
 f(t)=\sum_{{\nu}=0}^\infty f_{2\nu} R_{2\nu}(t)\in {\cal F}^A_{43}
 \end{equation}
is extremal iff inequality (\ref{h73}) turns into equality for this function. By the
property (\ref{h39}),  it is necessary for this that the coefficients $f_{\nu}$ in the
representation (\ref{f-extr}) of the extremal function $f$ be equal to zero for
\mbox{$\nu\ne 0,\ 1,\ 2,\ \ldots,\ 6,\ 9.$} In particular, this means that the function
$f$ is a polynomial of degree up to eighteen. Besides, the equality $L(f)=\lambda
(1)f(1)$ must hold. This equality is equivalent to the following conditions:
 $$
 f\left(\frac{1}{2}\right)=f(a_1)=f(a_2)=f(a_3)=0.
 $$
Thus, the point $\ds1/2$ is at least a zero and the points $a_1,\ a_2,\  a_3$ are double
zeros of the polynomial $f.$ Consequently, $f$ has the form
 \begin{equation}\label{f-extr7}
 f(t)=h(t)e(t),
 \end{equation}
 where
 $$
 h(t)= \left(t^2-\frac{1}{4}\right)\prod_{i=1}^{3}(t^2-a_i^2)^2,
 $$
and $e$  is a polynomial of degree up to four. Hence, by Lemma~\ref{lb3}, $f=cf^*,$ where
$c$ is a positive constant. This completes the proof of Theorem~\ref{tb1}.

\ \

 The author wishes to express his
thanks to V.V.~Arestov for suggesting the interesting problem, his attention to the work,
and help in the preparation of this paper. The author is also grateful to  A.G.~Babenko
for his active interest in the investigation.

\end{document}